\documentclass{amsart}
\usepackage{amsmath}
\usepackage{amsfonts}
\usepackage{latexsym}
\usepackage{amssymb}
\usepackage{enumerate}
\usepackage[dvips]{graphics}

\newcommand{\R}{{\mathbf R}}
\newcommand{\C}{{\mathbf C}}

\newcommand{\kk}{{\mathbf k}}
\newcommand{\cat}{{\rm {cat }}}

\newcommand{\id}{{\rm {id }}}

\newcommand{\tc}{{\rm\bf {TC}}}
\newcommand{\rk}{\rm {rk}}
\newcommand{\comment}[1]{}

\newcommand{\A}{{\mathcal A}}

\newcommand{\F}{{\mathcal F}}

\newcommand{\M}{{\mathcal M}}
\newcommand{\N}{{\mathcal N}}

\newcommand{\U}{{\mathcal U}}

\def\\N{\mathbb{N}}
\def\RR{\mathbf R}
\def\CC{\mathbf{C}}

\def\ZZ{\mathbb{Z}}
\def\R{\mathcal{FR}}

\def\R{\mathbf R}

\newtheorem{theorem}{Theorem}

\newtheorem{prop}[theorem]{Proposition}
\newtheorem{lemma}[theorem]{Lemma}

\newtheorem{corollary}{Corollary}
\newtheorem{definition}{Definition}

\numberwithin{equation}{section}

\begin{document}

\title[Collision Free Motion Planning]{Topological Robotics: Subspace Arrangements\\
and
Collision Free Motion Planning}

\author[M. Farber]{Michael Farber}
\address{Department of Mathematics, Tel Aviv University, Tel Aviv 69978, Israel}
\email{mfarber@tau.ac.il}
\thanks{M. Farber was partially supported by a grant from the US - Israel Binational Science Foundation
and by the Herman Minkowski Center for Geometry. S. Yuzvinsky is grateful to the Tel Aviv University for hospitality.}

\author[S. Yuzvinsky]{Sergey Yuzvinsky}
\address{Department of Mathematics
University of Oregon, Eugene OR 97403-1222, U.S.A. }
\email{yuz@math.uoregon.edu}


\date{October 8, 2002}

\dedicatory{To S.P. Novikov on the occasion of his 65-th birthday}


\begin{abstract}
We study an elementary problem of the topological robotics: collective 
motion of a set of $n$ distinct particles which one has to move 
from an initial configuration to a final configuration, with the requirement that 
no collisions occur in the process of motion.
The ultimate goal is to construct an algorithm which will perform this task once the initial and the final 
configurations are given. This reduces to a topological problem of finding the topological complexity
$\tc(C_n(\R^m))$ (as defined in \cite{F1, F2}) of the configutation space $C_n(\R^m)$ of $n$ distinct ordered particles 
in $\R^m$. We solve this problem for $m=2$ (the planar case) and for all odd $m$, including the case $m=3$ (particles in the three-dimensional space).
We also study a more general motion planning problem in Euclidean space with a hyperplane arrangement as obstacle.
\end{abstract}
\maketitle
\section{Introduction}

In this paper we study an elementary problem of the topological robotics: collective 
motion of a set of $n$ distinct particles which one has to move 
from an initial configuration to a final configuration, with the requirement that 
no collisions occur in the process of motion. This problem is clearly of practical interest. It becomes 
quite difficult when the number of 
particles is large. The ultimate goal is to construct an algorithm which will perform this task once the initial and the final 
configurations are given, see \cite{La}, \cite{Sh}. 

Any such motion planning algorithm must have instabilities \cite{F1}, i.e.
the motion of the system will always be discontinuous as a function of the initial and final configurations.
These instabilities of motion are caused by topological reasons.
A general approach to study instabilities of robot motion was suggested recently in \cite{F1, F2}.
With any path-connected topological 
space $X$ one associates a number $\tc(X)$, {\it the topological complexity of $X$}. 
This number is of fundamental importance for the 
motion planning problem: $\tc(X)$ determines the character of instabilities which have all  
motion planning algorithms in $X$. See \S \ref{section1} for a brief summary. 

In this paper we tackle the problem of calculating the topological complexity of the complements of subspace arrangements.
In particular, we compute the topological complexity of the 
configuration spaces of $n$ distinct points on the plane $\R^2$ and in the space $\R^3$. 
Our main results can be stated as follows: 

\begin{theorem}
\label{general1} Any motion planner for collision free motion of $n$
distinct points on the plane $\R^2$
has order of instability at least $2n-2$. There exist motion
planners having order of instability equal $2n-2$. 
\end{theorem}
\begin{theorem}
\label{general2} Any motion
planner for collision free motion of $n$ distinct points in the
three-dimensional space $\R^3$ has order of instability at least
$2n-1$. There exist motion planners having order of instability
equal $2n-1$.
\end{theorem}

\section{The motion planning problem}\label{section1}
\label{definition}

In this section we recall some definitions and results from \cite{F1, F2}. In particular we define the terms which are 
used in the statements of our main Theorems \ref{general1} and \ref{general2}.

Let $X$ be a topological space, thought of as the configuration
space of a mechanical system. 
Given two points $A, B \in X$,
one wants to connect them by a path in $X$. We always assume that $X$ is
a connected CW complex.
A solution to this motion
planning problem is a rule that takes $(A,B)
\in X \times X$ as an input and produces a path from $A$ to $B$ as an
output. Let $PX$ be the space of all continuous paths $\gamma: [0,1] \to X$,
equipped with the compact-open topology, and
let $\pi: PX \to X \times X$ be the map assigning the end points to a
path: $\pi (\gamma) = (\gamma(0),
\gamma(1))$. The map $\pi$ is a fibration whose fiber is the
based loop space $\Omega X$.
The motion planning problem consists of finding a section $s$ of this
fibration.

The section $s$ cannot be continuous, unless $X$ is
contractible, see \cite{F1}. One defines $\tc(X)$, the {\it
topological complexity of} $X$, as the smallest
number $k$ such that $X \times X$ can be
covered by $k$ open sets $U_1, \dots, U_k$, so that for every $i=1,
\dots, k$ there exists a continuous section $s_i:
U_i \to PX, \pi \circ s_i = \id$.

\begin{definition}\label{planner}
A motion planner in $X$ is given by finitely many subsets $F_1, \dots, F_k$ $\subset X\times X$ and by continuous maps
$s_i: F_i\to PX$, where $i=1, \dots, k$, such that:
\begin{enumerate}
\item[(a)] the sets $F_1, \dots, F_k$ are pairwise disjoint ($F_i\cap F_j=\emptyset$, $i\not= j$) and cover $X\times X$;
\item[(b)] $\pi\circ s_i =1_{F_i}$ for any $i=1, \dots, k$;
\item[(c)] each $F_i$ is an ENR.
\end{enumerate}
\end{definition}

The subsets $F_i$ are the {\it local domains} of the motion planner; the maps $s_i$ are the {\it local rules}.
Any motion planner determines a {\it motion planning algorithm}:
given a pair $(A, B)$ of initial - final configurations, 
we first determine the index $i\in \{1, 2, \dots, k\}$,
such that the local domain $F_i$ contains $(A, B)$; then we apply the local rule $s_i$ and 
produce the path $s_i(A,B)$ in $X$ as an output.

\begin{definition}
The order of instability of a motion planner is defined as the largest integer $r$ such that the closures of some $r$
among the local domains $F_1, \dots, F_k$ have a non-empty intersection: 
\begin{eqnarray*}
\bar F_{i_1}\cap \bar F_{i_2}\cap \dots \cap \bar F_{i_r}\, \not=\, \emptyset
\quad\mbox{where}\quad 1\leq i_1<i_2<\dots< i_r\leq k.
\end{eqnarray*} 
\end{definition}

The order of instability describes character of discontinuity of the motion planning algorithm
determined by the motion planner. 

\begin{theorem}\label{coincide} {\rm (\cite{F2})} Let $X$ be a connected smooth manifold. Then 
the minimal integer $k$, such that $X$ admits a motion planner
with $k$ local rules, 
equals $\tc(X)$.  Moreover,
the minimal integer $r>0$, such that $X$ admits a motion planner with order of instability $r$, equals $\tc(X)$.
\end{theorem}

This theorem explains importance of knowing the number $\tc(X)$ while solving practical motion planning problems.

Let us mention now some other results from \cite {F1,F2}, which will be used later in this paper.

\begin{theorem}\label{basic1}{\rm (\cite{F1, F2})} 
$\tc(X)$ depends only on the homotopy type of $X$. One has
$$
{\rm cat} (X) \leq \tc(X) \leq 2\ {\rm cat} (X) -1
$$
where ${\rm cat} (X)$ is the Lusternik-Schnirelmann 
category of $X$.
If $X$ is
\hfill\break
$r$-connected then
\begin{eqnarray}
\tc(X) < \frac{2\ {\rm dim}(X)+1}{r+1} + 1.\label{rconnected}
\end{eqnarray}
\end{theorem}

The following result provides a lower bound for $\tc(X)$ in terms of the
cohomology ring  $H^\ast(X;\kk)$ with coefficients in a
field $\kk$. The tensor product $H^\ast(X;\kk)\otimes H^\ast(X;\kk)$ is also a graded ring
with the multiplication
$$
(u_1\otimes v_1)\cdot (u_2\otimes v_2) = (-1)^{|v_1|\cdot |u_2|}\,
u_1u_2\otimes v_1v_2
$$
where $|v_1|$ and $|u_2|$ are the degrees of the cohomology classes
$v_1$ and $u_2$. The cohomology multiplication
$H^\ast(X;\kk)\otimes H^\ast(X;\kk)\to H^\ast(X;\kk)\label{prod}$ is a ring
homomorphism. Let  $I \subset H^\ast(X;\kk)\otimes
H^\ast(X;\kk)$ be the kernel of this homomorphism. The ideal $I$ is
called {\it the ideal of zero-divisors} of
$H^\ast(X;\kk)$. The {\it zero-divisors-cup-length} is the length of the
longest nontrivial product in the ideal of
zero-divisors.

\begin{theorem}\label{basic2}{\rm (\cite{F1})}
The topological complexity $\tc(X)$ is greater
than the zero-divisors-cup-length of
$H^\ast(X;\kk)$.
\end{theorem}

The topological complexity $\tc(X)$, as well as 
the Lusternik-Schnirelmann category $\cat(X)$, are particular cases of the notion of {\it Schwarz
genus} (also known as {\it sectional category}) of a fibration; it was introduced and 
thoroughly studied by A. Schwarz in \cite{Sch}.

\section{Hyperplane arrangements as obstacles: the main results}

In this section we study the topological complexity of complex
hyperplane arrangements complements. An important special case, which we mainly have in mind, 
is given by the configuration space of $n$ distinct points on the plane 
$$\{(z_1, z_2, \dots, z_n); z_i\in \CC, z_i\not=z_j\}.$$

Let $\mathcal A=\{H\}$ be a finite set of hyperplanes in an affine complex space
$\CC^{\ell}$. We will denote by $M(\A)$ the complement, i.e. $M(\A) = \CC^{\ell} - \cup_{H\in \A} H$. 
We will study the motion planning problem in $M(\A)$. From a different point of view, we may say that
we live in $\CC^\ell$ and the union of hyperplanes $\cup H$ represent our obstacles.

If $\bigcap_{H\in\A}H\not=\emptyset$ then $\A$ is called {\it central},
and up to change of coordinates the hyperplanes can be assumed linear.

Suppose that $\A$ is linear.
For each $H\in\A$ one can fix a linear functional $\alpha_H$ (unique up to a
non-zero multiplicative constant) such that $H=\{\alpha_H=0\}$. 
A set of hyperplanes $H_i\in \A$ is called {\it linear independent} if the corresponding functionals $\alpha_{H_i}$
are linearly independent. 
The rank of $\{\alpha_H\}$, i.e., the cardinality of
a maximal independent subset, is called the {\it rank} of $\A$ and denoted by
$\rk(\A)$. Clearly $\rk (\A) \leq \ell$ and the equality occurs
if and only if $\bigcap_HH=0$. 

If $\A$ is not central we define its rank as the
rank of a maximal central subarrangement of $\A$.

If $\A_i$ ($i=1,2$) are central arrangements in spaces $\CC^{\ell_i}$ respectively
then one can define their {\it product} as the arrangement $\A=\A_1\times\A_2$ in the space
$\CC^{\ell_1}\oplus \CC^{\ell_2}$ consisting of $H\oplus \CC^{\ell_2}$ for $H\in\A_1$
and $\CC^{\ell_1}\oplus H'$ for $H'\in\A_2$. It is easy to see that
$\rk(\A) = \rk(\A_1)+\rk(\A_2)$ and $M(\A)=M(\A_1)\times M(\A_2)$.

While dealing with the arrangement complements we will need the following
nontrivial result - see \cite{OT}, Section 5.2: {\it if $\A$ is an arbitrary
arrangement of rank $r$ then $M(\A)$ has homotopy type of a
simplicial complex of dimension $r$.}

One of the most interesting series of examples of central arrangements are the
complexifications of the sets of mirrors for Weyl groups. They are called
reflection arrangements. For instance the reflection arrangement of type $A_{m-1}$
is given in $\CC^m$ by the equations $z_i-z_j=0$ for all $1\leq i<j\leq m$; it has
rank $m-1$. The complement of this arrangement (i.e. of the union of its
hyperplanes) can be identified with the configuration space of all $m$-tuples
of distinct points on the plane $\R^2$ or $\C$. This space appears in Theorem \ref{general1} in the
Introduction.

Now we state the main theorems about the topological complexity of the arrangement complements:

\begin{theorem}
\label{upper}
 Let $M$ be the complement of a central complex
hyperplane arrangement of rank $r$. Then the topological
complexity satisfies $\tc(M) \leq 2r$.
\end{theorem}

This is slightly better than the upper bound $2r+1$ given by (1) mentioning the fact (see above) that
the complement $M(\A)$ has homotopy type of an $r$-dimensional complex, (2) using the homotopy invariance of $\tc(X)$,
see \cite{F1}, and (3) invoking Theorem \ref{basic1}. 

The above estimate can be improved in the case where the
arrangement is a product.

\begin{theorem}
\label{upper-product}
 If $\A = \A_1\times  \A_2 \times \dots \times
\A_k$, where for every $i=1, 2, \dots, k$, $\A_i$ is a central complex hyperplane arrangement,
then
$$\tc(M(\A)) \, \leq\, 2\sum_{i=1}^k r_i -k +1 = 2r -k +1.$$
Here $r_i$ denotes the rank of $\A_i$ and $r=\rk (\A)$.
\end{theorem}

The next theorem gives a necessary condition when the upper bound given by Theorem \ref{upper} is exact:

\begin{theorem}\label{lower} Let $\A$ be a central complex hyperplane
arrangement of rank $r$. Assume that there exist $2r-1$
hyperplanes $H_1, H_2, \dots, H_{2r-1}\in \A$ such that $H_1,
H_{2}, \dots, H_{r}$ are independent and for
any $j=1, 2, \dots, r$ the hyperplanes
$H_j, H_{r+1}, \dots, H_{2r-1}$ are independent. Then
$$\tc(M(\A)) =2r.$$
\end{theorem}

There are at least two large classes of arrangements satisfying
conditions of the previous theorem:

1. {\bf Generic arrangements of cardinality at least $2r-1$}. These
are the central arrangements whose any subset of cardinality $r$
is independent. Then any subset of cardinality $2r-1$ satisfies
the condition.

2. {\bf Reflection arrangements for reflection groups of types
$A_n, B_n$, and $D_n$.} Since every arrangement of type $B_n$ and
$D_n$ contains a subarrangement of type $A_{n-1}$, it is enough to show that the condition is satisfied
for the $A_{n-1}$ arrangement. Recall that the rank of $A_{n-1}$ is $r=n-1$. 
Denote by $H_{ij}$ the hyperplane given by
the equation $z_i-z_j=0$. Then the set of hyperplanes 
$\{H_{1j},H_{2k}; j=2,3,\ldots,n, k=3,4,\ldots,n\}$ satisfies the 
condition of Theorem \ref{lower}.

\begin{corollary}\label{main}
Let $C_n(\R^2)$ denote the configuration space of ordered sequences of $n$ distinct 
points on the plane. Then the topological complexity
of $C_n(\R^2)$ equals $2n-2$, i.e. $\tc(C_n(\R^2))=2n-2$. \qed
\end{corollary}

This Corollary follows from Theorem 8 and the above discussion.

Corollary \ref{main} combined with Theorem \ref{coincide} implies our main Theorem \ref{general1}.

\section{Proofs of Theorems \ref{upper} and \ref{upper-product}}

In order to prove Theorem \ref{upper} we need to use the well-known relations
between the complement of a central arrangement and that of the projectivization
 of this arrangement.
Let $\{H_1,H_2,\ldots,H_n\}$ be a central arrangement of hyperplanes in
$\CC^{\ell}$ of rank $r$.
 Since $M=M(\A)$ is invariant with
respect of the $\CC^*$-action on $\CC^{\ell}$ by multiplication,
we can consider the factor-space $M^*=M/\CC^*$.
This space is nothing but the complement in $\CC{\bf P}^{\ell-1}$ to the union of
the projectivizations of all $H_i$. The easiest (although noncanonical) way to
represent $M^*$ is to choose the coordinates so that $H_1=\{z_{\ell}=0\}$ and put
$z_{\ell}=1$ in the equations of $H_2,\ldots, H_n$. We obtain a not necessarily central arrangement of rank $\leq r-1$,
whose complement is $M^*$. Hence we see that $M^\ast$ has homotopy type of a cell complex of dimension
$\leq r-1$. Using homotopy invariance of $\tc(X)$ and the dimensional upper bound (see Theorems 3 and 4 from \cite{F1})
we obtain that $\tc(M^*)\leq 2r-1$.

A simple and well-known
observation is that the fibration $M\to M^*$ is trivial, i.e. $M$ is homeomorphic
to $M^*\times \CC^*$ (cf. \cite{OT}, Proposition 5.1).
Now we may apply the product
inequality (see \cite{F2}, Theorem 11) combined with the obvious fact $\tc(\C^*)=2$. We obtain
$\tc(M)\leq 2r-1+2-1=2r$. This proves Theorem \ref{upper}.

Now we want to prove Theorem \ref{upper-product}.
Suppose $\A$ is a central arrangement of rank $r$
 that can be represented as the product
$\A_1\times\A_2\times \cdots\times \A_k$ of central arrangements. Denote by $M_i$ the
complement of $\A_i$ and by $r_i$ its rank. Using again the product
inequality from \cite{F2} and Theorem \ref{upper} we have
$$\tc(M)=\tc(M_1\times M_2\times\cdots\times M_k)\leq \sum_{i=1}^k 2r_i-k+1
=2r-k+1.$$                   
Here $r$ is the rank of $\A$.
\qed

\section{Proof of Theorem \ref{lower}}

\subsection*{Monomials and flags of flats}

In this section we will use only the combinatorics of a hyperplane arrangement
$\A$ coded in its matroid. This means that the only information about $\A$ we need
it is what its subsets are independent. In fact we can forget about the
arrangement and consider an arbitrary simple matroid (e.g., see \cite{W}).

Let $\M$ be a simple matroid of rank $r$ on a set $S$ and $R$ a commutative ring.
Recall that the Orlik-Solomon algebra of $\M$ over $R$ is a graded $R$-algebra
$A=\oplus_{p=1}^rA_r$ that is the factor of the exterior algebra of the free
$R$-module with basis $\{e_s|s\in S\}$ over the ideal generated by $d(e_{s_1}\wedge
e_{s_2}\wedge\cdots\wedge e_{s_p})$ for all dependent subsets
$\{s_1,\ldots,s_p\}\subset S$. Here $d$ is the differential of the exterior
algebra of degree -1 satisfying the graded Leibniz condition and sending every
$e_s$ to 1. If $\M$ is matroid of an arrangement then $A$ is naturally isomorphic
to $H^*(M(\A);R)$ (\cite{OT}, Section 5.4).

Now we
recall relations between nonzero monomials in $A_r$, ordered bases
of $\M$, and maximal flags of flats of $\M$. Here a base of $\M$ is a maximal
independent subset of $S$ and a flat is a subset of $S$ closed with respect to
the dependence relation. If $\M$ is
a matroid of an arrangement, then flats are just the
intersections of some hyperplanes ordered opposite to inclusions.

Each (linearly) ordered subset $T$ of $S$ of cardinality $p$
defines a monomial $m(T)=\prod_{t\in T}e_t$ in $A_p$.
Here the product is taken in the
order on $T$. First, $m(T)\not=0$ if and only if
$T$ is independent in $\M$. On the other hand independent $T$ generates a flag
$\F(T)$ in the lattice of all the flats of $\M$ ordered by inclusion. If
$T=(t_1,t_2,\ldots,t_p)$ then
$$\F(T)=(<t_p>\subset<t_p,t_{p-1}>\subset\cdots\subset<t_p,t_{p-1},\ldots
,t_1>)$$
where $<U>$ is the flat generated by a subset $U$ of $S$. Notice
that the rank of the $i$-th flat in the flag is $i$.

The nonzero monomials in $A_p$ generate $A_p$ as a linear space
and are linearly dependent in general. A way to define a monomial
basis of $A_p$ is as follows. Fix an order on $S$ and call a
subset $C$ of $S$ a broken circuit if it is independent but there
is $s\in S- C$ such that $s<s'$ for every $s'\in C$ and
$C\cup\{s\}$ is a circuit (i.e., a minimal dependent set). Then a
subset $T\subset S$ is an {\bf nbc}-set if it does not contain any
broken circuit. Providing all subsets of $S$ with the induced
orders we obtain the set of monomials $\{m(T)\}$ in $A_p$ where
$T$ is running through all the {\bf nbc}-sets of cardinality $p$.
This set which we denote by ${\bf nbc}_p$ form a basis of $A_p$
(cf. \cite{OT}, sect. 3.1). We emphasize that this set depends on
the order in $S$.

Going back to a flag $F=(X_0=\emptyset\subset X_1\subset
X_2\subset\cdots\subset X_p)$ of flats one notices that an
arbitrary choice $s_{p-i}\in F_i=
 X_{i+1}- X_{i}$
($i=0,\ldots,p-1$) produces an ordered independent set
 $T=(s_1,\ldots,s_p)$ whence a nonzero
monomial $m(T)$. We denote the set of all these ordered bases by
$\nu(F)$ and the respective set of monomials by $\mu(F)$. Notice
that for each $T$ as above we have $\F(T)=F$. A monomial from
$\mu(F)$ is in ${\bf nbc}_p$ if and only if the chosen elements
are the smallest (in a fixed order on $S$) in each $F_i$ and their
order in $T$ coincides with the induced order. In particular for a
flag $F$ there is at most one ${\bf nbc}$-monomial in
$\mu(F)$. We call this monomial $m(F)$. If it does not exist then we put
$m(F)=0$. A flag $F$ for which $m(F)\not=0$ is called standard
(cf. \cite{SV}).

The flags  $F=(X_0\subset X_1\subset\cdots\subset X_p)$ of flats
with ${\rm rk}(X_i)=i$ provide a convenient language in order to
describe the decomposition of a monomial $m(T)$ into a linear
combination of the $\bf nbc$-monomials. We denote by
$\U=\oplus\U_p$ the linear space on the set of all flags as a basis
over a field $K$
graded by their length.

The following lemma is probably known to specialists. H.Terao
has informed the second author that he has known it since 1995 but we could not
find a proof in the literature.
\begin{lemma}
\label{decomposition} Let $T$ be an independent set in $\M$ taken
with the order induced by the fixed order on $S$. Then
$$m(T)=\sum({\rm sgn} \sigma)m(\F(\sigma T))$$
where $\sigma$ is running through the symmetric group $\Sigma_p$
and is acting on $T$ by permuting its elements.
\end{lemma}
\proof We will use a small part of the
techniques from \cite{OT}, Section 3.4. We again denote by $A=\oplus A_p$
the graded Orlik-Solomon algebra of $\M$ with coefficients in $R$.
Define the graded linear map $\phi: \U\to A$ via $\psi(F)=m(F)$
for every flag $F$ from $\U$. Also define the graded linear map
$\psi:A\to\U$ via $\psi(m(T)) =\sum({\rm sgn} \sigma)\F(\sigma
T)$ for every independent subset $T$ of $S$ taken with the induced order. The
fact that $\psi$ is well-defined is proved among other things in
\cite{OT}, Lemma 3.107.

Now we claim that $\psi$ is a section of $\phi$. Indeed let $T$ be
an independent subset of $S$ from ${\bf nbc}_p$. Then for every
$\sigma\in\Sigma_p$ the minimal element of each flat of the flag
$\F(\sigma T)$ is from $T$. This implies that $m(\F(\sigma
T))\not=0$ if and only if $\sigma$ is the identity permutation
whence $\phi\psi(m(T))=m(T)$. Since the monomials $m(T)$ for those
$T$ form a basis of $A$ the result follows.
\qed

Let $\M$ be an arbitrary matroid of rank $r$ on $S$ and $A$ its OS
algebra. Consider the graded algebra $\bar A=A\otimes A$ where the
tensor product is in the category of graded algebras over the
basic field. In particular, if $b\in A_k$ and $c\in A_{\ell}$ then
$(a\otimes b)(c\otimes d) =(-1)^{k\ell}ac\otimes bd$ for all
$a,d\in A$. For each standard generator $e_s$ ($s\in S$) of $A$
define the element of degree one of $\bar A$ via 
$$\bar
e_s=1\otimes e_s-e_s\otimes 1.$$ 
Our goal is to study the product
$$\pi=\prod_{s\in S}\bar e_s$$
taken in some order on $S$. If
$|S|>2r$ then clearly $\pi=0$. In the rest of this section we
assume that $|S|\leq 2r$.

The following lemma is straightforward.

\begin{lemma}
\label{sum} (i) Fix an order on $S$. Then
$$\pi=\sum_{(T,T')}(-1)^{|T|}{\rm sign}(\sigma) \cdot m(T)\otimes m(T'),$$
where $(T,T')$ runs through all the pairs of complementary independent sets
taken in the induced orders and $\sigma$ is the shuffle on $S$ preserving
orders inside $T$ and $T'$ and putting every element of $T'$ after
all elements of $T$.

(ii) If an order on $S$ changes via a permutation $\tau$ then
$\pi$ gets multiplied by ${\rm sign}\tau$. \qed
\end{lemma}

Clearly the linear space $A_p\otimes A_q$ is generated by simple
tensors $m_1\otimes m_2$ where $m_1$ and $m_2$ are nonzero
monomials from $A_p$ and $A_q$ respectively, 
and a basis in this space is formed by the
simple tensors where both monomials are ${\bf nbc}$.

The following lemma is an immediate corollary of Lemma
\ref{decomposition}.

\begin{lemma}
\label{complement} Fix an order on $S$ and let $\bar m=(m_1,m_2)$
be a pair of monomials both ${\bf nbc}$. Order respectively
the sets $T_1$ and $T_2$ of their elements and consider the flags
$F^i=\F(T_i)$ ($i=1,2$). Then any simple tensor of monomials from
Lemma \ref{sum} having in its decomposition $\bar m$ with a
nonzero coefficient $c$ has the form $m(T)\otimes m(T')$ where $T$
is complementary to $T'$, $T\in\mu(F^1)$, and $T'\in \mu(F^2)$.
Moreover $c=\pm 1$. \qed
\end{lemma}

\subsection*{Nonvanishing products in the tensor square.}

It is not hard to prove directly that for any matroid of rank $r$ the product
of any $2r$ elements $\bar e_s$ is 0. Later we will get this result for arrangements as a corollary
of others.
In this section we consider matroids of rank $r$ with the number
of elements less than $2r$. For these matroids the product
$$\pi=\prod_i\bar e_i$$
can be nonzero in $A\otimes A$.

\begin{prop}
\label{nonzero} Suppose that $S$ is the disjoint union of two
sets $T_1$ and $T_2$ with the following
properties: $T_1$ is independent and  $T_2\cup\{e_s\}$
is independent for every $s\in S$.
 Then $\pi\not=0$.
\end{prop}
\proof Put $|S|=n$ and $|T_2|=p$. Clearly $p<r$ and $n-p\leq r$.
Consider the flag $F=\F(T_1)=(X_0=\emptyset\subset
X_1\subset\cdots\subset X_p)$ corresponding to some order on $S$.
Without any loss of generality we can assume that $F$ is standard
in this order and put $m=m(F)=m(T_2)$. The condition on $T_2$
implies that $|X_{i+1}- X_i|=1$ for all
$i=0,1,2,\ldots,p-1$. Thus $m(T_2)$ is the only monomial having
$m$ in its decomposition with a nonzero coefficient. If $m'$ is an
arbitrary ${\bf nbc}_{n-p}$-monomial from the decomposition of
$m(T_2)$ then $m(T_1)\otimes m(T_2)$ is the only simple tensor of
complementary monomials having $m\otimes m'$ is its decomposition
which completes the proof.              \qed

{\bf Proof of Theorem \ref{lower}.}
Let $\A$ be an arrangement satisfying the condition of Theorem \ref{lower}.
Denote by $S$ the subarrangement consisting of $H_1,\ldots,H_{2r-1}$ and by $\M$
the matroid on $S$ of this subarrangement. By Proposition \ref{nonzero},
$\pi=\prod_{s\in S}\bar e_s\not=0$. Since the Orlik-Solomon algebra $A'$ of any
subarrangement is a natural subalgebra of the Orlik - Solomon algebra $A$ of the whole
arrangement (see \cite{OT}, Proposition 3.66), the product $\prod_{s\in S}\bar e_s\not=0$ in the cohomology ring of
$M(\A)$. Applying the cohomological lower bound for the topological complexity (see Theorem \ref{basic2}), we obtain
$\tc(M(\A))\geq 2r$. This and Theorem \ref{upper} imply the result.     \qed

The proofs of the following corollaries are straightforward.

\begin{corollary}
\label{affine}
Let $\A$ be an arrangement satisfying the condition of Theorem \ref{lower}.
Then $\tc(M(\A)^*)=2r-1$. \qed
\end{corollary}

\begin{corollary}
\label{vanishing}
For an arbitrary complex central hyperplane arrangement $\A$ of rank $r$
the product of any $2r$ elements from the kernel of the multiplication homomorphism 
$H^*(M)\otimes H^*(M)\to H^*(M)$ equals zero. \qed
\end{corollary}

\section{Example: motion planner for collision free control of three particles}
\label{3points}

Here we will describe a recipe to obtain an explicit motion planner for moving a triple of ordered
points in $\CC$ (or $\RR^2$) avoiding collisions. 

The configuration space of such
triples 
$$C_3= \{(z_1, z_2, z_3); z_1\not=z_2, z_1\not=z_3, z_2\not=z_3\}$$ 
has superfluous coordinates: the first point $z_1$ may be arbitrary 
and one has to observe only the relative positions of the second and the third particles with respect to the first. 
Hence $C_3$ is homeomorphic (via the map $(z_1, z_2, z_3) \mapsto (z_1, z_2-z_1, z_3-z_1)$)
to the product $\C\times M$, 
where 
$$M=\CC^2- (H_1\cup H_2\cup H_3),$$
with $\CC^2=\{(z_1,z_2),\ z_1, z_2\in\CC\}$, $H_1=\{z_1=0\}$, $H_2=\{z_2=0\}$, and $H_3=\{z_1=z_2\}$.

One may perform another change of coordinates and represent $M$ as a product $M^\ast\times \C^\ast$,
where $M^\ast = \C-\{0,1\}$ and $\C^\ast=\C-\{0\}$. 
We assign to a pair $(z_1, z_2)\in M$,  
$$h(z_1, z_2) = \displaystyle{\left(\frac{z_1}{z_2}, z_2\right)}\in M^\ast\times \C^\ast.$$ 
The obtained map $h: M\to M^\ast \times \C^\ast$ is a homeomorphism. 

Now we see that $M^\ast$ is homotopy equivalent to a one-dimensional complex, hence it admits a motion planner 
with 3 sets, see Theorems 3 and 4 in \cite{F1}. 
A motion planner on $M^\ast$ can be easily described explicitly.
$\C^\ast=\C-\{0\}$ admits a motion planner with two sets, see \cite{F1}; such motion planner may also be explicitly
constructed. The explicit construction of a motion planner on the product $M^\ast\times \C^\ast$ may be now 
obtained by repeating the 
arguments used in the proof of the product inequality (see Theorem 11 of \cite{F1},  or the construction after Theorem 12.1 in \cite{F2}). 
This gives a motion planner in $M$ having 4 local domains. Making this recipe precise is straightforward.

\section{Motion planning for collision free motion of $n$ particles in the 3-dimensional space}

In this section we study the problem of constructing a 
collision free control for $n$ distinct particles in $\RR^3$.
Our main result is:

\begin{theorem}
\label{R}
Let $M=C_n(\R^3)$ be the configuration space of $n$ distinct ordered particles in $\RR^3$.
Then
$\tc(M)= 2n-1$.
\end{theorem}

\begin{proof}
Similarly to the plane
case, the configuration space of $n$ ordered distinct points in $\RR^3$
can be viewed as the complement to an arrangement of linear subspaces. Indeed the
space of all $n$-tuples of points from $\RR^3$ can be identified with $\RR^n\otimes
\RR^3=\RR^{3n}$ with the coordinates $x_{ik}$, where $i=1,2,\ldots, n$ and $k=1,2,3$. Here
the numbers $x_{j1}, x_{j2}, x_{j3}$ are the 
space coordinates of the particle number $j$.
To exclude coincidences one needs to remove from $\R^{3n}$ the union of the subspaces
$H_{ij}$ ($1\leq i<j\leq n$), where $H_{ij}$ is given by the system of 3
equations 
$$x_{i1}-x_{j1}=0, \quad x_{i2}-x_{j2}=0, \quad x_{i3}-x_{j3}=0.$$ 
Hence the configuration
space we are interested in is $M=\RR^{3n}-\bigcup H_{ij}$.

It follows immediately from this representation of $M$ that it is simply connected. Indeed, $M$ is obtained from the Euclidean space
by removing finitely many subspaces of codimension 3. 

To find a low and upper bounds for the topological complexity
$\tc(M)$ it is
convenient to start with its cohomology ring $H^*(M)$ with the integral coefficients. For a general subspace
arrangement this ring is much more complicated than Orlik-Solomon algebras and is not
even defined by the combinatorics of the arrangement. The property of our arrangement
that simplifies the matter is the fact that the dimensions of all our subspaces and
their intersections are divisible by 3.  Moreover the intersection
lattice of the arrangement in $\RR^{3n}$ coincides with the lattice of the braid
arrangement corresponding to $n$ distinct points in $\CC$. Under these conditions
the cohomology algebra $A=H^*(M)$ is defined similarly to the Orlik-Solomon algebra with the only difference that the degree of the natural generators
$e_{ij}$ is $2$, see \cite{Yu}, Section 7. This implies that $A$ is a commutative algebra. More precisely
$A=\ZZ[e_{ij}]/I$, where the ideal $I$ is generated by the polynomials $e_{ij}^2$ and
$$e_{ij}e_{ik}-e_{ij}e_{jk}+e_{ik}e_{jk}$$
for every triple $1\leq i<j<k\leq n$.

Except being commutative, $A$ has properties very similar to the respective
Orlik-Solomon algebra. For instance, $A=\oplus_{p=0}^{n-1}A_p$ although here
degree of $A_p$ is $2p$. A monomial of the form $e_{i_1j_1}e_{i_2j_2}e_{i_kj_k}$ is nonzero in $A$ if and only if
the respective linear functionals $x_{i_\alpha}-x_{j_\alpha},$ where $\alpha = 1, \dots, k$,
are linearly independent.

In particular, we find that the top dimension, where $M$ has a nontrivial cohomology, is $2(n-1)$.
Since $M$ is simply connected, it follows that $M$ is homotopy equivalent to a cell-complex of dimension 
$\leq 2(n-1)$ (see for example, \cite{EG}).
Then using inequality (\ref{rconnected}) we obtain
$$\tc(M)<\frac{{4(n-1)+1}}{2}+1,$$
i.e. $\tc(M)\leq 2n-1$. 

Now we will use the cohomological lower bound for the topological complexity given by Theorem
\ref{basic2}.
We will find a non-zero
product in $A\otimes A$ of $2n-2$ elements of the form
$$\bar e_{ij}= 1\otimes e_{ij}-e_{ij}\otimes 1,$$
which are zero-divisors. 
Consider the following product 
$$\pi=\prod_{i=2}^n\, (\bar e_{1i})^2\, \in A\otimes A.$$
An easy computation gives $(\bar e_{ij})^2=-2e_{ij}\otimes e_{ij}$ for arbitrary $i,j$.
Hence we find 
$\pi=(-2)^{n-1}m\otimes m$, where 
$$m=\prod_{i=2}^ne_{1i}.$$ 
Since the
linear functionals $\{x_1-x_i|i=2,3,\ldots,n\}$ are linear independent, 
the monomial $m\not=0$ is nonzero in $A$, and hence the product $\pi$ is nonzero.
Thus we obtain the opposite inequality $\tc(M)\geq 2n-1$. \qed
\end{proof}

Theorem \ref{general2} follows from Theorem \ref{R} combined with Theorem \ref{coincide}.

{\bf Remark.} Let $C_n(\R^m)$ denote the configuration space of $n$ ordered distinct points in $\R^m$.
Repeating the argumant of Theorem \ref{R} we find $\tc(C_n(\R^m))=2n-1$ for $m$ odd.

{\bf Conjecture.} {\it For $m$ even $\tc(C_n(\R^m)) = 2n-2$}. 

This would generalize our Corollary \ref{main} where $m=2$. Our arguments show that for any even $m$ 
the topological complexity $\tc(C_n(\R^m))$ equals either $2n-2$ or $2n-1$.

\end{document}